\documentclass[journal]{IEEEtran}
\usepackage{graphicx,epsf,psfrag,amsmath,amssymb,amsfonts,verbatim}
\usepackage[mathscr]{eucal}
\usepackage{algorithmic}
\usepackage{verbatim,pifont,color,epsf,psfrag,amssymb,amsfonts,amsmath}
\usepackage{caption}
\usepackage{subcaption}
\usepackage{accents}
\IEEEoverridecommandlockouts                              

\newtheorem{theorem}{Theorem}

\def\nn{\nonumber}
\def\beq{\begin{equation}}
\def\eeq{\end{equation}}
\def\bea{\begin{eqnarray}}
\def\eea{\end{eqnarray}}
\def\ba{\begin{array}}
\def\ea{\end{array}}
\def\defeq{{\stackrel{\Delta}{=}}}

\def\eg{{\it e.g., \/}}

\def\ie{{\it i.e.,\ \/}}

\definecolor{bgrd}{rgb}{1,1,1}
\definecolor{grey}{rgb}{0.9,0.9,0.6}
\definecolor{gray}{rgb}{0.5,0.5,0.5}

\makeatletter
\newdimen{\captionwidth}
\long\def\@makecaption#1#2{%
\captionwidth .9\hsize
\vskip 10pt%
\setbox\@tempboxa\hbox{#1: #2}%
  \ifdim \wd\@tempboxa >\captionwidth%
    \setbox\@tempboxa\hbox{#1:\hspace*{.5em}}%
    \hfil\parbox{\captionwidth}{\raggedright\hangindent \wd\@tempboxa%
    \hangafter=1\unhbox\@tempboxa#2}\hfill%
  \else\centerline{\box\@tempboxa}%
  \fi
}
\makeatother

\def\T{\mbox{\small T}}

\def\span{\mbox{span}}
\def\tr{\mbox{trace}}

\newcommand{\mbbE}{\mathbb{E}}

\newcommand{\Pmsc}{\mathscr{P}}

\def\figwidth{.5\textwidth}

\def\T{{\mbox{\tiny T}}}
\def\u{{\mbox{\tiny u}}}
\def\d{{\mbox{\tiny d}}}
\begin{document}
%
\title{Dynamic Pricing and Distributed Energy Management for Demand Response}
%
%
%

\author{Liyan Jia and Lang Tong, {\em Fellow, IEEE}

\thanks{Liyan Jia and Lang Tong are with the School of Electrical and Computer Engineering, Cornell University, Ithaca, NY, USA {\tt\small \{lj92,lt35\}@cornell.edu}

This manuscript has been accepted by IEEE Transactions on Smart Grid.

This work is supported in part by the National Science Foundation under Grant CNS-1135844 and Grant 15499. %
}%

}

%



\maketitle


\begin{abstract}
  The problem of dynamic pricing of electricity in a retail market is considered. A  Stackelberg game is used to model interactions between a retailer and its customers; the retailer sets the day-ahead hourly price of electricity and consumers adjust real-time consumptions to maximize individual consumer surplus.

  For thermostatic demands, the optimal aggregated demand is shown to be an affine function of the day-ahead hourly price. A complete characterization of the trade-offs between consumer surplus  and retail profit is obtained.  The Pareto front of achievable trade-offs is shown to be concave, and each point on the Pareto front is achieved by an optimal day-ahead hourly price.

  Effects of integrating renewables and local storage are analyzed. It is shown that  benefits of renewable integration all go to the retailer when the capacity of renewable is relatively small. As the capacity increases beyond a certain threshold, the benefit from renewable that goes to consumers increases.  
\end{abstract}


\section{Introduction}\label{sec:intro}
Demand response in a smart grid is expected to offer economic benefits to consumers while improving overall operation efficiency and reliability. A properly designed demand response program can  reduce the peak load, compensate for uncertainties associated intermittent renewables, and reduce the cost of system operation.

Two types of demand response are commonly used.  One gives the retail utility direct control of consumers' consumptions.  See, \eg \cite{DemandResponsePJM,Chao&Etal88,CLRC:05}.  Although this form of demand response provides flexibilities to the system operator, a consumer loses the ability to manage consumption based on her own preferences.

An alternative is to reshape demand response through dynamic pricing. Examples of this type of programs include the use of  critical pricing \cite{Faruqui&etal:09} and real time pricing \cite{Borenstein&etc:02,Borenstein&Holland:05}.  Such schemes allow the consumers to mange consumptions individually.

In this paper, we consider a demand response scheme of the second type.    We focus on a particular type of dynamic pricing referred to as the {\em day-ahead hourly pricing} (DAHP).  DAHP was first considered in \cite{Borenstein&etc:02,Borenstein&Holland:05} and has already been implemented by some utility companies in U.S. \cite{Hopper&Goldman&Neenan:05EJ}. Under DAHP, the hourly retail prices of electricity are set one day ahead of the actual consumptions, thus providing  a level of price certainty to consumers.  DAHP also allows  the retailer to adjust prices on a day-to-day basis, taking into account operating conditions at the wholesale market. The empirical study of DAHP reported in \cite{Hopper&Goldman&Neenan:05EJ} concludes that DAHP ``not only improves the linkage between wholesale and retail markets, but also promotes the development of retail competition.''

\vspace{-1em}
\subsection{Summary of Main Results}
We analyze in this paper interactions between a retailer and its customers.  We use the term
retailer to include the traditional retail utilities as well as energy aggregators in a smart grid setting.
  A Stackelberg game model is used to capture retailer-consumer interactions in which  the retailer is  the leader who sets the  DAHP and the consumers the followers who adjust individual consumption.  For thermostatically controlled load,  we show in Section~\ref{ssec:demand} that the optimal aggregated demand is an affine function of the DAHP, and the optimal DAHP can be obtained via convex optimizations.

\begin{figure}[tb]
\begin{center}
\begin{psfrags}
\psfrag{cs}[c]{$\textsf{\text{\footnotesize Consumer Surplus (CS)}}$}
\psfrag{rp}[l]{$\textsf{\text{\footnotesize Retail Profit (RP)}}$}
\psfrag{xo}[l]{$\pi^{\text{o}}$}
\psfrag{xsw}[l]{$\pi^{\text{sw}}$}
\psfrag{xr}[l]{$\pi^{\text{r}}$}
\psfrag{d}[r]{$\Delta$}
\psfrag{0}[l]{$0$}
\includegraphics[width=2.7in]{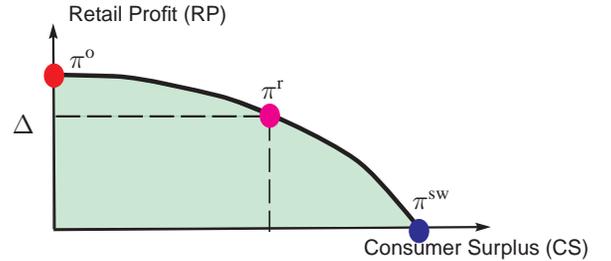}
\end{psfrags}
\caption{\small CS-RP trade-off curve with different dynamic pricing schemes.}
\label{fig:trade-off1}
\end{center}
\end{figure}

Main results of this paper are built upon the  characterization of the Pareto front of trade-offs between consumer surplus (CS) and retail profit (RP), as illustrated in Fig.~\ref{fig:trade-off1}, where each point on the Pareto front is associated with an optimized DAHP. We show in Section~\ref{sec:Game} that the Pareto front is concave and monotonically decreasing.  Several well known pricing schemes are shown in Fig.~\ref{fig:trade-off1}. In particular, in the absence of renewable integration and energy storage, the  price $\pi^{\text{sw}}$  that maximizes the sum of consumer surplus and retail profit (the social welfare) is  shown to result in zero retail profit.  This means that the social welfare maximization price is not economically viable to the retailer. The optimal regulated monopoly price $\pi^{\text{r}}$ is associated with the point on the Pareto front where the retailer's profit is fixed at some $\Delta$. The retail profit maximization price is $\pi^{\text{o}}$ when no constraint is imposed on the retailer. We note also that benchmark pricing schemes such as constant pricing, time of use (ToU) and proportional markup pricing are suboptimal; they are in general inside the Pareto front, indicating that improvement can be achieved by applying the optimization framework developed in Section~\ref{sec:tradeoff}.

We also consider effects of incorporating renewable energy and energy storage.  We show that, when the retailer integrates renewable sources, the retail profit at the social welfare optimizing DAHP is positive, making social welfare maximization potentially a viable operation objective. It is shown that the benefit of renewable integration all goes to the retailer when the capacity of renewable integration is small. As the capacity increases, so does the benefit to consumers. For the consumer side storage, we show that the CS-RP trade-off can be formulated under the same Stackelberg game framework.

Simulation studies are presented using realistic electricity prices and home energy management models to illustrate the benefits of distributed energy management and DAHP-based demand response. The effects of renewable generation on the trade-off curve between consumer surplus and retail profits are illustrated via simulations.

\vspace{-0.5em}
\subsection{Related work}
 The literature on dynamic pricing at retail level appeared along with the beginning of wholesale market deregulation. In \cite{Borenstein&etc:02, Borenstein:05}, the authors demonstrated the benefit of DAHP and showed that DAHP could attract consumers in the long run.  In \cite{Carrion&Etal:07TPS, Conejo&Etal:08TPS}, retail price optimization was formulated by considering consumer side uncertainties and  financial risks. Relative to this line of existing work,  the main contribution of this paper is a full characterization of the trade-off between the achievable retail profit and consumer surplus.

The use of Stackelberg game model to study interactions   between a retailer and its consumers goes back to  \cite{Luh&etal:82TAC}. The authors of \cite{Luh&etal:82TAC} presented a general game theoretic framework and a load adaptive pricing scheme.  The pricing scheme considered in  \cite{Luh&etal:82TAC} is different from the DAHP treated here, primarily because the load adaptive pricing was designed for a vertically integrated utility whereas the DAHP considered in this paper is applicable for a deregulated market operating based on a two-settlement mechanism.    Instead of aiming at the single objective of global social optimality as in \cite{Luh&etal:82TAC}, we give a full characterization of optimal solutions for a set of (local) objectives.


Extensive results exist on optimizing demand response operations.   For example, a distributed algorithm was developed  to achieve social welfare maximization in \cite{Li&Etal:11PESGM}. A hierarchical market structure and distributed optimization is also considered in \cite{Joo&Ilic:13TSG}.  Contract design between renewable generators and the aggregator is studied in \cite{Papavasiliou&Oren:11PESGM}, where the aggregator is responsible for a large scale PHEV charging. A particularly relevant work, from a modeling perspective, is the generalized battery model \cite{Hao&Etal} that captures succinctly the thermostatic control as a generalized battery model. The demand model used in this paper falls in this category.

\vspace{-0.5em}
\subsection{Organization and notations}
This paper is organized as follows. In Section \ref{sec:Game}, a Stacelberg game is introduced and rational behaviors of the retailer and its customers  are analyzed. Section \ref{sec:tradeoff} characterizes the tradeoff between the retail profit and consumer surplus. In Section \ref{sec:renewable},  renewable integrations by the retail utility is considered.  Section~\ref{sec:storage} considers the problem when storage devices are used at the consumer side.  In Section \ref{sec:sim}, we conduct numerical simulations to compare different benchmark schemes and investigate the effects of renewable integration.

 Throughout the paper, $(x_1,\cdots, x_N)$ is for a column vector and $[x_1,\cdots, x_M]$ a row vector. We use $x^+$ to represent $\max\{0,x\}$, $x^-$ to represent $\max\{0,-x\}$, and $\bar{x}$ the expected value $\mathbb{E}[x]$.


\section{A Stackelberg game model}
\label{sec:Game}

\subsection{DAHP and a Stackelberg game model}
 We model the retailer-consumer interactions through DAHP using a Stackelberg game where the retailer is the leader by setting the DAHP and  consumers the follower by adjusting their individual consumption.

\begin{enumerate}
\item
The  DAHP $\pi=(\pi_1,\cdots, \pi_{24})$ is a 24-dimensional real vector where $\pi_i$ is the retail price of electricity at hour $i$. The DAHP is set one day ahead and fixed throughout the day of operation. The pay-off function of the retailer will be defined later in Section~\ref{ssec:retail}.
\item
In real time, a consumer dynamically determines her consumption. If storage is available, the consumer may also charge or discharge the storage. Selling back surplus power in a net metering setting is allowed. The payment from a consumer to the retailer is settled as the inner product of $\pi$ and the real-time energy consumption. The pay-off function of consumers is the consumer surplus that measures the satisfaction of the consumer. One candidate is provided in Section~\ref{ssec:demand}.

\item
The retailer meets the aggregated demand by purchasing electricity at the wholesale market, possibly complementing the purchase by its own renewables.
\end{enumerate}


The Stackelberg game can be solved via a backward induction.  We thus present first the analysis of consumers' response to a fixed DAHP $\pi$. The problem of optimizing $\pi$ from the retailer's point of view is then considered.


\subsection{Consumer action: optimal demand response}
\label{ssec:demand}
We start our analysis of the Stackelberg game from the consumer side. After the DAHP $\pi$ is given, a consumer dynamically determines her own consumption to optimize the consumer surplus defined by the difference between the utility of consumption and the electricity payment. It is assumed that the consumer in a smart grid can adjust her consumption based on real-time measurements. The optimal control policy and the resulted demand response serve as the predicted behavior of a rational consumer in the Stackelberg game described above.

In this paper, we consider price elastic demands that are results of optimal control of a linear system.  A typical example is the thermostatically controlled loads.  See \cite{Hao&Etal} for a general model for other types of applications.

For thermostatic control involving heating-ventilation and air conditioning (HVAC) units, empirical study \cite{Bargiotas&Birddwell:88ITPD} shows that the dynamic equation that governs the temperature evolution is given by
\begin{equation}
\begin{array}{r l}
x_i&=x_{i-1}+ \alpha(a_i-x_{i-1}) - \beta p_i + w_i, \\
y_i& = (x_i, a_i) + v_i,
\label{eq:hvacmodel}
\end{array}
\end{equation}
where $i$ is the index of period, $y_i$ the noisy measurements of the indoor temperature $x_i$ and outdoor temperature $a_i$, $p_i$ the power drawn by the HVAC unit, $w_i$ the process noise and $v_i$ the measurement noise. Both $w_i$ and $v_i$ are assumed to be zero mean and Gaussian. System parameters $\alpha$ ($0<\alpha<1$) and $\beta$ model the insolation of the building and efficiency of the HVAC unit. Note that the above equation applies to both heating and cooling scenarios. But we exclude the scenario that the HVAC does both heating and cooling during the same day. We focus herein on the cooling scenario ($\beta > 0$) and the results apply to heating ($\beta < 0$) as well.

In the following, we will use superscript $(j)$ to denote the variables associated with consumer $j$. For consumer $j$, assume that she wants to keep the indoor temperature close to her own desired temperature $t_i^{(j)}$ for the $i$th hour. The deviation of the actual indoor temperature $x_i^{(j)}$ from the desired temperature $t_i^{(j)}$ can be used to measure consumer $j$'s discomfort level. By assuming symmetric upward and downward discomfort measure, a quadratic form of consumer utility over $N$ periods  is given by
\begin{equation}
u^{(j)}= -\mu^{(j)} \sum_{i=1}^{N}  (x_i^{(j)}-t_i^{(j)})^2,
\label{eq:utility}
\end{equation}
where $\mu^{(j)}$ is consumer $j$'s weight factor to convert the temperature deviation to a monetary value.

The consumer surplus of consumer $j$ is defined by the difference between the consumer utility and the total payment from the consumers. In particular, given a DAHP $\pi=(\pi_1,\cdots, \pi_{24})$ and energy consumption $p_i^{(j)}$ in the $i$th hour, the payment for the real-time consumption from consumer $j$ is settled as $\sum_{i=1}^{24} \pi_i p_i^{(j)}$.  Therefore, her consumer surplus  is a random variable  given by

\begin{equation}
\textsf{\text{cs}}^{(j)} \defeq u^{(j)} - \pi^{\T} p^{(j)},
\label{eq:cs}
\end{equation}
where $p^{(j)}=(p^{(j)}_1,\cdots, p^{(j)}_{24})$.

We use the consumer surplus as the payoff function of consumers in the Stackelberg game, and the action of a consumer is to decide the electricity consumption during each period. As an assumption, a rational consumer maximizes consumer surplus in real time according to DAHP and her current energy state $x_i^{(j)}$. Therefore, the solution to the following stochastic control program characterizes the demand response to DAHP.

\begin{equation}
\begin{array}{r l}
\max_{p^{(j)}} & \mathbb{E}\left\{\sum_{i=1}^{24} [- \mu^{(j)} (x_i^{(j)} - t_{i}^{(j)})^2 ] - \pi^{\T} p^{(j)} \right\} \\
\mbox{s.t.} & x^{(j)}_{i} = x^{(j)}_{i-1} + \alpha^{(j)} (a_{i} - x^{(j)}_{i-1}) - \beta^{(j)} p^{(j)}_i + w^{(j)}_i, \\
& y^{(j)}_i=(x^{(j)}_i,a_i)+ v^{(j)}_i . \\
\end{array}
\label{eq:opt_d}
\end{equation}


Under mild conditions where the price doesn't vary too much during a day and $\mu^{(j)}$ is large, we ignore the positivity constraint and the rate constraint for energy consumption $p^{(j)}$. The backward induction gives the  optimal demand response as shown in the following theorem first appeared in \cite{Jia&Tong:12Allerton}.
\begin{theorem}
\label{thm:opt_demand}
Under DAHP, the optimal aggregated residential demand response for a fixed DAHP $\pi$ has the following form,
\begin{equation}
d(\pi)= -G \pi+b,
\label{eq:demand}
\end{equation}
where $d(\pi)$ is a random vector whose $i$th entry is the aggregated demand from all consumers in hour $i$, $G$ a {\em deterministic} and positive definite matrix, depending only on system parameters and user preferences, and $b$ a Gaussian random vector corresponding to the overall process and observation noise.
\end{theorem}

{\em Proof: }  see the Appendix. \hfill $\Box$

\vspace{0.5em}
Theorem~(\ref{thm:opt_demand}) establishes an affine relationship between the optimal demand response and day-ahead hour price under some mild conditions. It is shown that the sensitivity matrix of demand with respect to price, $-G$, is not affected by the realization of randomness. This means that the change of price will have deterministic effect on the expected value of demand.
Also, we  see that the retailer does not need to estimate each consumer's parameters in Eq.~(\ref{eq:opt_d}), including $\mu$, $\alpha$, $\beta$ and etc.. It only needs to fit the affine model as in Eq.~(\ref{eq:demand}) based on historical data.  Herein,  we assume that rational consumers follow the optimal demand response according to Eq.~(\ref{eq:demand}).

Given the retail price $\pi$, by substituting the optimal demand response in Theorem~\ref{thm:opt_demand} back into the consumer optimization problem, the expected total consumer surplus  (summed over all consumers)  is given by
\bea
\overline{\textsf{\text{cs}}}(\pi) & =&  \sum_j \overline{\textsf{\text{cs}}}^{(j)}(\pi), \nn\\
&=& \sum_j \mbbE\left([- \sum_{i=1}^{24} \mu (x^{*(j)}_i - t_{i})^2] - \pi^{\T}p^{*(j)}(\pi),\right) \nn\\
&=& \pi^{\T}G\pi/2 -\pi^{\T}\bar{b} + c, \label{eq:avg_cs}
\eea
where we used the fact that $\mbbE(p^{*(j)}(\pi))=-G^{(j)}\pi+\bar{b}^{(j)}$ for the $j$th consumer and
$G \defeq \sum_j G^{(j)}$,  $\bar{b}\defeq \sum_j \bar{b}^{(j)}$.  Here
$x^{*(j)}$ and $p^{*(j)}(\pi)$ are the same as the values in the proof of Theorem~\ref{thm:opt_demand} in the Appendix, and $c$ is a constant depending on the variance of the noise.

\subsection{Retailer's action: optimal dynamic pricing}
\label{ssec:retail}
In this paper, we assume that the retailer is a price taker in the wholesale real time market.  This means that the aggregated demand from consumers in real time does not affect the wholesale price. Additionally, we assume that the Stackelberg game discussed in this paper is one with perfect information, which means the leader has complete knowledge of the follower's payoff function. In the real-time operation, the retailer is required to satisfy aggregated demand by paying for the distribution cost and the cost of procuring power from the wholesale market.

Let $\lambda = (\lambda_1, \lambda_2,...,\lambda_{24})$ denote the random vector of average marginal cost during each hour, which includes the wholesale electricity price, distribution cost and other service cost.  The total expected daily retail profit $\textsf{\text{rp}}(\pi)$, defined as the difference between the expected real time retail revenue and the retail cost, is given by
\bea
\overline{\textsf{\text{rp}}}(\pi) &=& \mbbE\left((\pi-\lambda)^{\T}d(\pi),\right) \nn\\
&=& (\pi - \bar{\lambda})^{\T}(-G \pi+\bar{b}), \label{eq:avg_rp}
\eea
where we use the assumption that, given fixed DADP $\pi$,  the real-time wholesale price is independent of the real-time consumption.

The retailer's pricing decision depends on its own payoff function. A particularly relevant objective is the social welfare defined as the sum of consumer surplus and retail profit, $\ie$
\begin{equation}
\overline{\textsf{\text{sw}}}(\pi) = \overline{\textsf{\text{rp}}}(\pi) + \overline{\textsf{\text{cs}}}(\pi).
\end{equation}
Using (\ref{eq:avg_cs}-\ref{eq:avg_rp}), we obtain the following theorem that characterizes the social welfare maximizing DAHP.
\begin{theorem}
\label{thm:SW}
The optimal retail price $\pi^{\text{sw}}$ that maximizes the expected social welfare is the expected real time retail cost,
\[
\pi^{\text{sw}}=\bar{\lambda},
\]
and the expected retail profit under $\pi^{\text{sw}}$ is
$\overline{\textsf{\text{rp}}}(\pi^{\text{sw}})=0$.  And
for any $\pi'$ such that $\overline{\textsf{\text{rp}}}(\pi') \ge 0$, $\overline{\textsf{\text{cs}}}(\pi') \le \overline{\textsf{\text{cs}}}(\pi^{\text{sw}})$.
\end{theorem}

{\em Proof:} See the Appendix. \hfill $\Box$

\vspace{0.5em}
Theorem~\ref{thm:SW} shows that, if the social welfare is to be maximized, the retailer generates no profit.   It is also shown that when social welfare is the payoff function, the retailer simply matches the DAHP with the expected retail cost.

\section{Pareto Front of Tradeoffs}
\label{sec:tradeoff}
In this section, we solve the Stackelberg game with a weighted social welfare as the retailer's payoff function.  Our goal is to characterize the Pareto front of the CS vs. RP tradeoffs.   In particular, we consider the following optimization for the retailer
\begin{equation}
  \label{eq:optII}
  \mbox{max}~ \{\overline{\textsf{\text{rp}}}(\pi)+ \eta \overline{\textsf{\text{cs}}}(\pi)\},
\end{equation}
where $\eta$ is a parameter that allows the retailer to weigh its profit against consumer's satisfaction.

 In practice,  $\eta$  can be set according to the operating cost, long term business plan, and social impact etc. The two extremes are
 the social welfare optimization when $\eta=1$ and profit maximization when $\eta=0$. For a profit regulated monopoly, the retailer sets $\eta$ to maintain the allowed level of profit. In this case, the retailer may want to maximize consumer surplus subject to obtaining its regulated profit in order to maintain or attract additional customers.

For the defined Stackelberg game, the action of the retailer is to determine the day-ahead hourly price to maximize the pay-off function. Given the expected cost $\bar{\lambda}$ of procuring electricity from the wholesale market, and the preference parameter $\eta$,  we can solve Eq.~(\ref{eq:optII}) by substituting (\ref{eq:avg_cs}-\ref{eq:avg_rp}) in (\ref{eq:optII}), yielding
\bea
\overline{\textsf{\text{rp}}}(\pi)+ \eta \overline{\textsf{\text{cs}}}(\pi)
&=& (\frac{\eta}{2}-1)\pi^{\T}G\pi \nn\\
& & +\pi^{\T}((1-\eta)\bar{b}+G\bar{\lambda}) + \eta c,\nn
\eea
from which we obtain the optimal price as
\begin{equation}
\label{eq:opt_pi}
\pi^{*}(\bar{\lambda}, \eta) = \frac{1}{2-\eta}\bar{\lambda} + \frac{1-\eta}{2-\eta}G^{-1}\bar{b}.
\end{equation}
Note that the optimized DADP in (\ref{eq:opt_pi}) is made of two terms. The first depends only on the expected wholesale price of electricity, which represents the cost to the retail. As an economic signal to consumers, this term corresponds to the behavior that the retail price follows the average wholesale price of electricity. The second term depends only on the consumer preference and randomness associated with consumer environments.

With (\ref{eq:opt_pi}), the Pareto front of CS-RP tradeoffs is characterized by
\beq \label{eq:P}
\Pmsc_{\bar{\lambda}}=\Bigg\{ \Big(\overline{\textsf{\text{cs}}}(\pi^{*}(\bar{\lambda},\eta)),\overline{\textsf{\text{rp}}}(\pi^{*}(\bar{\lambda},\eta))\Big): \eta \in [0,1]\Bigg\},
\eeq
where each element of $\Pmsc_{\bar{\lambda}}$ can be computed in closed form.

\vspace{0.5em}
The shape of the trade-off region is characterized by the following theorem.
\begin{theorem}
\label{thm:property}
The  Pareto front of $(\overline{\textsf{\text{cs}}},\overline{\textsf{\text{rp}}})$ is concave and decreasing. The area above the Pareto front is infeasible for the retailer to achieve under DAHP whereas the closed area on and below the Pareto front can be achieved by DAHP.
\end{theorem}

{\em Proof:} See the Appendix.  \hfill $\Box$

\vspace{0.5em}
We can now revisit Fig.~\ref{fig:trade-off1} in lights of the analysis developed in this section.  As summarized in Sec~\ref{sec:intro}.A, the achievable tradeoff region is convex in the CS-RP plane.  The Pareto front $\Pmsc_{\bar{\lambda}}$ is a collection of CS-RP tradeoffs achieved by optimal DADP. The social welfare maximizing pricing $\pi^{\text{sw}}$ results in zero retail profit, and profit regulated  monopoly price $\pi^{\text{r}}$ is located at the Pareto front where the retail profit has a regulated profit  $\Delta$.  Prices such as constant prices, ToU, and markups on day-ahead wholesale price are typically inside the Pareto front $\Pmsc_{\bar{\lambda}}$, thus suboptimal in general.  See validations from numerical simulations in Sec~\ref{sec:sim}.

An alternative formulation is based on  a practical situation where the retailer optimizes its profit under the constraint that the consumer surplus exceeds a certain level. In particular, the problem is formulated as
\begin{equation}
\label{eq:optI}
\max_\pi \overline{\textsf{\text{rp}}}(\pi)~~\mbox{subject to}~~\overline{\textsf{\text{cs}}}(\pi) \ge \tau.
\end{equation}
The following theorem shows the equivalence between (\ref{eq:optI}) and (\ref{eq:optII})  in obtaining the Pareto front in (\ref{eq:P}).

\begin{theorem}
  \label{thm:equivalence}
  For any specific $\eta$, if the solution in (\ref{eq:optII}) is $\pi^*$, $\pi^*$ is also a solution to (\ref{eq:optI}) with $\tau = \overline{\textsf{\text{cs}}}(\pi^*)$.  Varying $\tau$ in optimization (\ref{eq:optI}) and varying $\eta$ in optimization (\ref{eq:optII}) results in the same Pareto front $\Pmsc_{\bar{\lambda}}$ in (\ref{eq:P}).
  \end{theorem}

{\em Proof:} See the Appendix.  \hfill $\Box$

\section{Effects of renewable integration}
\label{sec:renewable}
We consider in this section the role of renewable energy at the retailer side. As a large load aggregator, the retailer may have the financial ability and incentive to have its own or  contracted solar/wind farms. We restrict ourself to the scenario that the retailer owned renewables are used to compensate the real time loads; the extra renewable power is spilled.

Denote the vector of hourly marginal cost of renweable as $\nu$. It is reasonable to assume  $\lambda - \nu > 0$ to hold every where in practice. Let $q = ( q_1,...,q_{24})$ be the maximum renewable power accessible to the retailer in each hour. Accordingly, the retailer's profit with renewable integration is changed to
\begin{equation}
\label{eq:opt_w}
\begin{array}{rcl}
\textsf{\text{rp}}_{\mbox{\tiny{Renew}}}(\pi) & = & \pi^{\T}d(\pi) - \min_{0 \le \tilde{q} \le q } \{ (\lambda- \nu)^{\T}(d(\pi) - \tilde{q})^+ \\
&&  + \nu^{\T} \tilde{q}\},\\
& = & \pi^{\T}d(\pi) - \nu^{\T}d(\pi)-  (\lambda - \nu )^{\T}(d(\pi) - q)^+,
\end{array}
\end{equation}
where $\tilde{q}$ is the actual renewable power the retailer used, and the function $(x)^+$ is the positive part of $x$, defined as $\mbox{max}\{x,0\}$.

Notice that the expected retail profit is also a concave function of the DAHP vector, $\pi$. Therefore, if the retailer's objective is profit maximization, the optimal price can be easily solved by a convex program.

Following similar discussions in Section \ref{sec:tradeoff}, we focus on the problem that the retailer's payoff function is weighted social welfare. Notice here that the consumer surplus does not change while the retail profit $\textsf{\text{rp}}(\pi)$ is replaced by $\textsf{\text{rp}}_{\mbox{\tiny{Renew}}}(\pi)$. Solving the following optimization problem will give the Pareto front for the RP vs. CS trade-offs in the presence of renewable integration.

\begin{equation}
\label{eq:opt_ew}
\mbox{max}~ \{\overline{\textsf{\text{rp}}}_{\mbox{\tiny{Renew}}}(\pi)+ \eta \overline{\textsf{\text{cs}}}(\pi)\}.
\end{equation}

Intuitively the achievable trade-off region between CS and RP will be enlarged due to renewables toward upright in the CS-RP plane. The following theorem verifies this intuition and further shows how the benefit of renewable integration is distributed between the retailer and its consumers.
\begin{theorem}
Assume that for each hour $i$, the renewable power is uniformly distributed with maximum capacity $K$, $\ie$ $q_i \sim U[0, K]$. For each preference parameter $\eta$, let the optimal price in (\ref{eq:opt_ew}) be $\pi_{\mbox{\tiny{Renew}}}^{\eta}$ and the corresponding expected demand $\bar{d}(\eta)$.
Define $\Delta \overline{\textsf{\text{rp}}}(\eta)$ and $\Delta \overline{\textsf{\text{cs}}}(\eta)$ as the increases of retail profit and consumer surplus due to renewable integration, respectively. We then have the following statements:

\begin{enumerate}
  \item Renewable integration always increases retail profit, $\ie$ $\overline{\textsf{\text{rp}}}_w(\pi_w^{\eta}) > 0$ for all $\eta$.
  \item When $K \le \min_{i \in \{1,...,24\}} \bar{d}_i(\eta)$, $\Delta \overline{\textsf{\text{cs}}}(\eta) = 0$. Otherwise, $\Delta \overline{\textsf{\text{cs}}}(\eta) > 0$.
  \item As $K \rightarrow \infty$, the fraction of renewable integration benefit to consumer surplus $\frac{\Delta \overline{\textsf{\text{cs}}}(\eta)}{\Delta \overline{\textsf{\text{cs}}}(\eta)+ \Delta \overline{\textsf{\text{rp}}}(\eta)} \rightarrow \frac{1}{3-2\eta}$.
\end{enumerate}
\label{thm:wind}
\end{theorem}

{\em Proof:} See the Appendix. \hfill $\Box$

\vspace{0.5em}
Theorem~\ref{thm:wind} shows that the benefit of renewable all goes to the retailer when the capacity is small. The intuition is that the retailer would naturally use the low cost energy to fulfill the need of demand and pocket the benefit.
As the capacity of renewable integration increases, the additional renewable will be spilled unless the consumption increases.  The only way to increase the retail profit is to increase consumption by lowering the price of electricity.  As a result, all consumer benefits from the price reduction. As $K \rightarrow \infty$, the fraction of renewable integration benefit to CS converges to a particular limit, $\frac{1}{3-2\eta}$, depending only on the weighting factor $\eta$.

The above theorem also shows  that when the retailer cares only about its own profit, the fraction of renewable integration to CS converges to $\frac{1}{3}$. On the other hand, if the retailer's objective is social welfare maximization, as the capacity of renewable power goes to infinity, the fraction of renewable integration benefit to consumer converges to 1 while the retail profit converges to zero.

\section{Effects of storage at the consumer side}

\label{sec:storage}

In this section, we consider the role of storage at the consumer side. In particular, we assume the net-metering option where a consumer can sell back the excess energy with the same purchasing price. In the presence of storage, a consumer can arbitrage over the hourly varying DAHP.

For hour $i$, denote the energy level in the battery as $B_i$ and the energy charged to the battery as $r_i$ (when $r_i \le 0$, it means discharging from the battery). Let $r_i = r_i^+ - r_i^-$ where $r_i^+ \ge 0$ and $r_i^- \ge 0$ are the positive and negative parts of $r_i$, respectively.   Then the dynamics of the battery can be expressed as
\begin{equation}
B_{i+1} = \kappa(B_{i} + \tau r_i^+ - r_i^-/\rho),
\end{equation}
where $\kappa \in (0,1)$ is the storage efficiency, $\tau \in (0,1)$ the charging efficiency and $\rho \in (0,1)$ the discharging efficiency. These efficiency factors model the potential energy loss during the process of storage, charging and discharging.

After incorporating the storage devices, the optimal demand response problem at the consumer side is changed to
\begin{equation}
\begin{array}{r l}
\underaccent{p,r,B}{\max} & \mathbb{E}\left\{\sum_{i=1}^{24} [- \mu (x_i - t_{i})^2] - \pi^{\T}(p + r) \right\} \\
\mbox{s.t.} & x_{i} = x_{i-1} + \alpha (a_{i} - x_{i-1}) - \beta p_i + w_i, \\
& y_i=(x_i,a_i)+ v_i, \\
& B_{i+1} = \kappa(B_{i} + \tau r_i^+ - r_i^-/\rho), \\
& B_{24} = B_{0}, 0 \le B_{i} \le C, \\
& 0 \le r_i^+ \le r^{\u}, 0 \le r_i^- \le r^{\d},
\end{array}
\end{equation}
where $B_0$ is the initial energy level in the storage, $C$ the capacity of the battery, $r^{\u}$ the charging limit, and  $r^{\d}$ the discharging limit.

Under the net-metering assumption, the optimization problem can be divided into two independent sub problems, where the first one is the same as the previous optimal stochastic HVAC control (\ref{eq:opt_d}) and the second one is purely energy arbitrage\cite{Xu&Tong:14PESGM}. This means that adding storage on the demand side doesn't change the original linear relationship between the actual HVAC consumption and retail price; the benefit of storage goes to the consumers only in the form of arbitrage options.

Therefore, given day-ahead price $\pi$, the optimal charging vector $r^*(\pi)$ can be solved via the following deterministic linear program:
\begin{equation}
\begin{array}{r l}
\underaccent{r, B}{\max} & -\pi^{\T}r  \\
\mbox{subject to} & B_{i+1} = \kappa(B_{i} + \tau r_i^+ - r_i^-/\rho), \\
& B_{24} = B_{0}, 0 \le B_{i} \le C, \\
& 0 \le r_i^+ \le r^{\u}, 0 \le r_i^- \le r^{d}.
\end{array}
\end{equation}
Using the original form of RP and CS without storage, $\textsf{\text{rp}}(\pi)$ and $\textsf{\text{cs}}(\pi)$, the retailer's payoff function in (\ref{eq:optII}) changes to
\begin{equation}
\label{eq:opt_b}
\mbox{max}~~\{\overline{\textsf{\text{rp}}}(\pi)+(\pi - \lambda)^{\T}r(\pi) + \eta ( \overline{\textsf{\text{cs}}}(\pi)- \pi^{\T} r(\pi))\}.
\end{equation}
Solving (\ref{eq:opt_b}) will give the induced tradeoff curve between CS and RP under different weighting factor $\eta$.

\section{Numerical Simulations}
\label{sec:sim}
In this section, we present simulation results that help to gain insights into effects of demand response with optimized DAHP. The parameters used in the simulations were extracted from  actual temperature record in Hartford, CT, from July 1st, 2012 to July 30th, 2012.  For the same period, we used the record of day-ahead wholesale price for the same location from ISO New England. The parameters for a HVAC thermal dynamic model (\ref{eq:hvacmodel}) were set at $\alpha = 0.5$, $\beta = 0.1$, $\mu=0.5$. The desired indoor temperature was set to be $18^{\circ}C$ for all hours. We assumed zero marginal renewable power cost.

\subsection{Benchmark comparisons}
\label{ssec:benchmark}

We first present direct comparisons between the optimized DAHP with some simple benchmarks. Specifically, we considered the following  well known schemes:

\begin{itemize}
  \item Constant pricing (CP): in this case, the price remained constant for the whole day, $\ie$ $\pi_1=\pi_2=...=\pi_{24}=x$.  Each value of $\pi$ had a corresponding CS and RP pair.  By varying $x$, we traced the performance of constant pricing on the CS-RP plane.
  \item Time of Usage (ToU): in this case, a single day was divided into two parts: peak hours and normal hours. For a normal hour $i$, $\pi_i = \pi^{\mbox{\tiny norm}}$, and for a peak hour $j$, $\pi_j = \pi^{\mbox{\tiny peak}}$, where $\pi^{\mbox{\tiny peak}} > \pi^{\mbox{\tiny norm}}$. In this section, we fixed $\pi^{\mbox{\tiny peak}} = 1.2 \pi^{\mbox{\tiny norm}}$ and set the peak hours as 9am to 5pm.

  \item Proportional mark-up pricing (PMP): in this case, the retail pricing was indexed by the day-ahead wholesale price.  Specifically, the ratio of DAHP at each hour over the day-ahead wholesale price at the same hour was constant, $\ie$ $\frac{\pi_1}{\mathbb{E}\lambda_1}=\frac{\pi_2}{\mathbb{E}\lambda_2}=...=\frac{\pi_{24}}{\mathbb{E}\lambda_{24}}=\gamma$. By varying $\gamma$, we traced the performance of this pricing on the CS-RP plane.
\end{itemize}

Qualitatively, we expected that the CS-RP performance of these schemes will fall inside the region defined by the Pareto front.

\begin{center}
\begin{figure}[!ht]
\includegraphics[width=\figwidth]{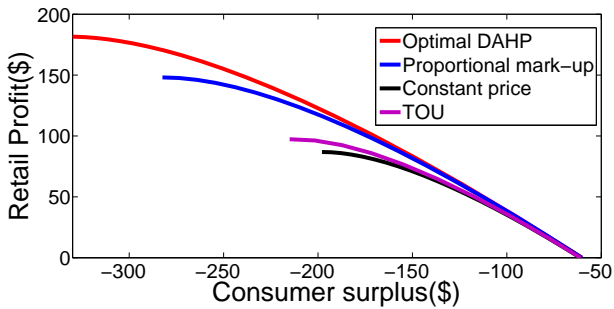}
\caption{Comparison of three pricing schemes}
\label{fig:cs_rp_tradeoff_nowind}
\end{figure}
\end{center}

The trade-off curves under the four schemes were plotted in Fig.~\ref{fig:cs_rp_tradeoff_nowind}. Theoretical statements established in the paper were evident:  the trade-off curves between CS and RP were downward and concave. The social welfare optimal operating point ($\eta=1$)  for the optimal DAHP resulted in zero retail profit.  For other pricing schemes, the points with zero profit were quite close to the social welfare optimal operating point under DAHP. Especially, since the social welfare maximization price was equal to the proportional price with $\gamma=1$, the social welfare optimal operating point was also on the trade-off curve under PMP. As expected, the traces of CP, PMP and {\tr ToU} all fell under the optimal DAHP which defines the Pareto front.

\begin{table}
\begin{center}
\begin{tabular}{|r|c|c|c|c|}
\hline
Regulated RP & DAHP & PMP & ToU & CP \\
\hline
85 & -152.0 & -154.3 & -170.2 & -183.4 \\
\hline
0 & -60.3 & -60.3 & -60.9 & -61.1 \\
\hline
\end{tabular}
\label{tab:compare}
\caption{Comsumer surplus under different regulated profits with four pricing schemes}
\end{center}
\end{table}

Table I illustrates quantitatively the gain of the optimized DAHP over typical benchmarks. The retail profits were fixed and associated consumer surplus values were read from the CS-RP tradeoff curves for all the pricing schemes. The performance of PMP was close to DAHP since PMP also took advantage of varying environment during a day, whereas TOU and CP performed much worse due to the lack of flexibility; DAHP achieved $10.7\%$ gain over ToU and $17.1\%$ gain over CP when the regulated RP was fixed as 85.

We also compared the payoffs to consumers with and without optimal demand response.  Here we assumed a baseline control policy aimed at maintaining the indoor temperature below the desired temperature plus some tolerance. Table~II shows the results of optimal demand response and the scenarios with tolerances set to be 0 and 2 degrees, under the social welfare maximization price. From the results we can see that with the increase of tolerance, the payment to retailer decreased but the uncomfort level grew fast. The optimal demand response was the best if measured by the consumer surplus as defined in Eq. (\ref{eq:cs}).


\begin{table}
\begin{center}
\begin{tabular}{|r|c|c|c|}
\hline
 & DR & Tolerance = 0 & Tolerance = $2^{\circ}C$ \\
\hline
Payment & 37.7 & 40.8 & 28.0 \\
\hline
Discomfort Level & 1.6 & 0 & 33.2 \\
\hline
Consumer Surplus & -39.3 & -40.8 & -61.2 \\
\hline
\end{tabular}
\label{tab:optimalresponse}
\caption{Comparison of consumer payoff with and without optimal demand response}
\end{center}
\end{table}

\subsection{Characteristics of DAHP}  
\begin{center}
\begin{figure}[!ht]
\begin{psfrags}
\psfrag{eta}[c]{$\eta$}
\includegraphics[width=\figwidth]{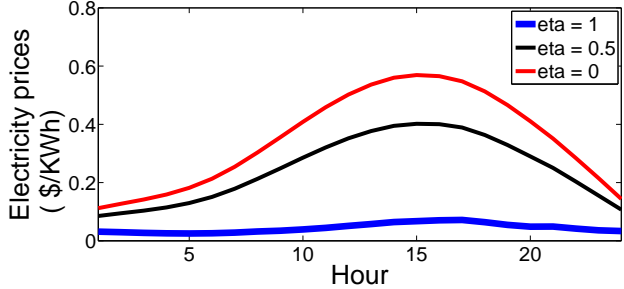}
\end{psfrags}
\caption{Comparison of consumer payoff with and without optimal demand response}
\label{fig:prices}
\end{figure}
\end{center}

To gain insights into how the optimal DAHP balances the tradeoff between RP and CS, we plotted the DAHP at each hour for different values of $\eta$ in Fig.~\ref{fig:prices}.  The two extreme points were $\eta=0$ for profit maximization and $\eta=1$ for social welfare optimization. In fact, the DAHP curve with $\eta=1$ was exactly the expected retail cost. We can see that, the RP maximizing DAHP showed a temporal pattern that matched the demand-supply dynamics.  The prices were higher at peak hours.   What was not shown in this figure was that the  consumers used less energy in response to higher prices and consumer surplus decreased. The CS-RP pair moved northwest along the Pareto front. In contrast, the social welfare optimal price showed a more consistent pricing across all hours.

\subsection{Renewable integration}
\label{ssec:sim_w}
\begin{center}
\begin{figure}[h!]
    \includegraphics[width =\figwidth]{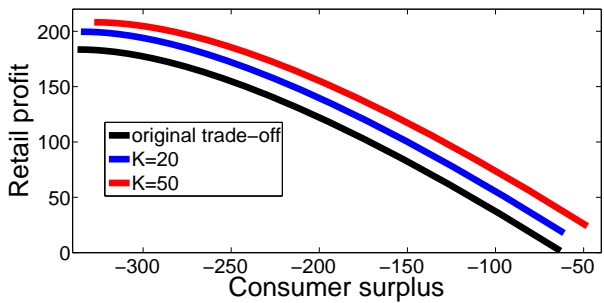}
    \caption{Trade-off curve with renewable integration}
    \label{fig:cs_rp_tradeoff_wind}
\end{figure}
\end{center}

Assuming renewable power in each hour was uniformly distributed over $[0,K]$, by solving the problem (\ref{eq:opt_ew}), we plotted the trade-off curves between CS and RP in Fig.~\ref{fig:cs_rp_tradeoff_wind}. Not surprisingly, renewable integration enlarged the tradeoff curve and all parties benefited from the low cost energy. In particular, the social welfare optimal pricing became economically viable, $\ie$ the social welfare maximization prices (rightmost points on the trade-off curve) resulted in positive retail profit. Furthermore, when the capacity of renewable was small ($K=20$), the trade-off curve moved upward, indicating that almost all the benefit from renewable integration went to the retailer. When the capacity was larger ($K=50$), the trade-off curve went upright, which showed that some part of the renewable integration benefit went to consumers.

\begin{center}
\begin{figure}[h!]
  \begin{psfrags}
\psfrag{eta}[c]{$\eta$}
\psfrag{ratio}[c]{$\frac{\Delta \overline{\textsf{\text{cs}}}(\eta)}{\Delta \overline{\textsf{\text{cs}}}(\eta)+ \Delta \overline{\textsf{\text{rp}}}(\eta)}$}
    \includegraphics[width = \figwidth]{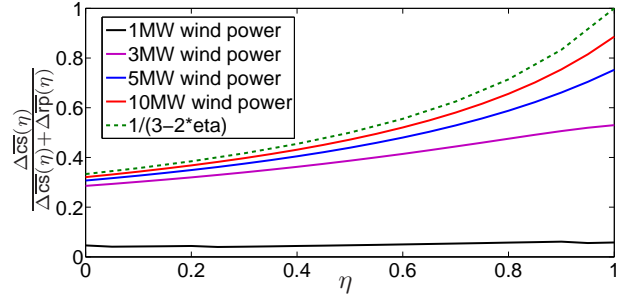}
    \end{psfrags}
    \caption{Fraction of renewable benefit to consumer}
    \label{fig:wind_benefit_dist}
\end{figure}
\end{center}

Fig.~\ref{fig:wind_benefit_dist} provided more insights into how the benefits of renewable were distributed among the retailer and the consumers.  As the level of integration ($K$) increased, so did the fraction of the renewable integration benefit to CS, $\frac{\Delta \text{cs}(\eta)}{\Delta \text{rp}(\eta)}$,  and it converged to $\frac{1}{3-2\eta}$ as $K$ goes to infinity.  More interestingly, as $\eta$ increased, more emphasis was placed on consumer surplus, substantial gains were achieved by the consumers.

\section{Conclusion}

In this paper, we have studied a day-ahead hourly pricing (DAHP) mechanism for distributed demand response in uncertain and dynamic environments. Such a pricing scheme has the advantage of reducing consumer anxiety of pricing uncertainties and allowing the retail utility to optimize the retail pricing adaptively.  The main contribution is a full characterization of the tradeoff between consumer surplus and retail profit for a class of demands.  This result allows the retailer to optimize its pricing strategy, taking into consideration his access to renewable energy.

A number of important issues not treated in this paper have been or are being addressed.  For example, the system parameters in the affine mapping of the aggregated load  are assumed to be known here.  In practice, obtaining these parameters require machine learning techniques.  See \cite{Jia&Tong&Zhao:13Allerton}.  The model considered here does not take into account the  closed-loop behavior in the sense that demand response in the distribution system has an impact in the wholesale market.  This is a challenging topic deserves further study. The developed theory can also be used  to study the effects of demand side storage and renewable integration. See \cite{Jia&Tong:15PESGM}.



\bibliographystyle{IEEEtran}
%


\bibliographystyle{IEEEbib}
{
\bibliography{Journal,Conf,Misc,Book}
}

\section*{APPENDIX}

\subsection*{Proof of Theorem~\ref{thm:opt_demand}}
For simplicity, in this paper we consider the case that control period length is 1 hour. For general cases where control period length is different from 1 hour, similar proof can be followed with some transformations. See \cite{Jia&Tong:13CDC} for details.

Through backward induction, the optimal control for consumer $j$  can be obtained as
\begin{equation*}
\begin{tabular}{r c l}
 $p_i^{*(j)}(\pi)$  &$=$&$ \frac{1}{\beta}\left(\hat{x}_{i-1|i-1}^{(j)}+\alpha^{(j)}(\hat{a}_{i|i-1}-\hat{x}^{(j)}_{i-1|i-1})-x_i^{*(j)} \right) $,\\
$x_i^{*(j)}$ & $\defeq$ & $\frac{\pi_i-(1-\alpha^{(j)})\pi_{i+1}}{2\mu^{(j)} \beta^{(j)}} + t_i^{(j)}$,
\end{tabular}
\label{eq:solution}
\end{equation*}
where $\hat{x}_{i-1|i-1}^{(j)}$ and $\hat{a}_{i|i-1}$ are the minimum mean squared error estimate of indoor and outdoor temperatures of hour $i$ at hour $i-1$, respectively.  Here $\pi_{25}$ is assumed to be 0, and $x_i^{*(j)}$ is an ancillary value. Expanding the solution above, the total demand of consumer $j$ is
\[
p^{*(j)}(\pi) = -G^{(j)}\pi + b^{(j)},
\]
where $b^{(j)}$ is Gaussian, independent of $\pi$, and
\[
  G^{(j)}_{ik} = \left\{
  \begin{array}{l l}

    [1+(1-\alpha^{(j)})^2]/[2\mu^{(j)} (\beta^{(j)})^2], & \quad \text{if $i=k\neq1$};\\
    1 / [2\mu^{(j)} (\beta^{(j)})^2], & \quad \text{if $i=k=1$};\\
    (-1+\alpha^{(j)})/[2\mu^{(j)} (\beta^{(j)})^2], & \quad \text{if $|i-k|=1$};\\
    0, & \quad \text{o.w.}\\
  \end{array} \right.
\]

Notice that for $0 < \alpha < 1$,
\[
1 / [2\mu^{(j)} (\beta^{(j)})^2] > |(-1+\alpha^{(j)})/[2\mu^{(j)} (\beta^{(j)})^2]|,
\]
\[
[1+(1-\alpha^{(j)})^2]/[2\mu^{(j)} (\beta^{(j)})^2] > 2|(-1+\alpha^{(j)})/[2\mu^{(j)} (\beta^{(j)})^2]|.
\]
Therefore, $G^{(j)}$ is deterministic and diagonal dominant with positive diagonal elements. Hence, $G^{(j)}$ is positive definite.

On the other hand, the optimal aggregated demand
\[
d(\pi)=\sum_k p^{*(j)}(\pi)=-G\pi+b,
\]
where $b =\sum_j b^{(j)}$, $G = \sum_j G^{(j)}$. $G$ is positive definite and deterministic, depending only on the parameters.

\subsection*{Proof of Theorem~\ref{thm:SW}}
Setting the derivative of $\overline{\textsf{\text{sw}}}(\pi)$  to zero gives the optimal price and resulted retail profit as,$\pi^{\text{sw}} = \bar{\lambda}$, and $ \overline{\textsf{\text{rp}}}(\pi^{\text{sw}}) = 0$. For any $\pi'$ such that $\overline{\textsf{\text{rp}}}(\pi') \ge 0$, we have
\[
\overline{\textsf{\text{cs}}}(\pi') = \overline{\textsf{\text{sw}}}(\pi') - \overline{\textsf{\text{rp}}}(\pi') \le \overline{\textsf{\text{sw}}}(\pi^{\text{sw}}) - 0 = \overline{\textsf{\text{cs}}}(\pi^{\text{sw}}).
\]

\subsection*{Proof of Theorem~\ref{thm:property}}
For $\eta \in [0,1]$, the solution to (\ref{eq:optII}) is given by
\[
\pi^* = \frac{1}{2-\eta}G^{-1}[(1-\eta)\bar{b} + G\bar{\lambda}].
\]
Define the resulted retail profit and consumer surplus as $\textsf{\text{rp}}^*(\eta) \defeq \overline{\textsf{\text{rp}}}(\pi^*(\eta))$, $\textsf{\text{cs}}^*(\eta) \defeq \overline{\textsf{\text{cs}}}(\pi^*(\eta))$. Then,
\[
\frac{\partial \textsf{\text{rp}}^*(\eta)}{\partial \textsf{\text{cs}}^*(\eta)} = \frac{\frac{\partial \textsf{\text{rp}}^*(\eta)}{\partial \eta}}{\frac{\partial \textsf{\text{cs}}^*(\eta)}{\partial \eta}} = -\eta.
\]
Because $\textsf{\text{cs}}^*(\eta)$ is an increasing function of $\eta$, $\frac{\partial \textsf{\text{rp}}^*(\eta)}{\partial \textsf{\text{cs}}^*(\eta)}$ decreases as $\textsf{\text{cs}}^*(\eta)$ increases. The curve is concave. According to Theorem~{\ref{thm:equivalence}}, $\textsf{\text{rp}}^*(\eta)$ is the optimal value of (\ref{eq:optI}) when consumer surplus is at least $\textsf{\text{cs}}^*(\eta)$. Therefore, no CS-RP pair can be above the trade-off curve.

\subsection*{Proof of Theorem~\ref{thm:equivalence}}
With a particular $\eta$, assume $\pi^*$ is the solution to (\ref{eq:optII}). Let $\tau = \overline{\textsf{\text{cs}}}(\pi^*)$ in (\ref{eq:optI}). Then $\pi^*$ will be in the feasible set of (\ref{eq:optII}). If there exists $\pi'$, such that $\overline{\textsf{\text{rp}}}(\pi') > \overline{\textsf{\text{rp}}}(\pi^*)$, and $\overline{\textsf{\text{cs}}}(\pi') \ge \tau$,
\[
\overline{\textsf{\text{rp}}}(\pi')+ \eta \overline{\textsf{\text{cs}}}(\pi')> \overline{\textsf{\text{rp}}}(\pi^*) + \eta \tau = \overline{\textsf{\text{rp}}}(\pi^*)+ \eta \overline{\textsf{\text{cs}}}(\pi^*).
\]
Hence, $\pi^*$ is not the solution to (\ref{eq:optII}) since $\pi'$ achieves better objective value. It contradicts with the assumption. Therefore, $\pi^*$ is also a solution to (\ref{eq:optI}).

\subsection*{Proof of Theorem~\ref{thm:wind}}
Before renewable integration, for any $\eta$, the first order condition gives that the optimal demand level $d(\eta)$ satisfies
\[
b - (2-\eta)d(\eta)= G\lambda.
\]
After renewable integration, for a particular $\eta$, the first order condition gives that the optimal demand level $d_w(\eta)$ satisfies
\[
b - G\nu - (2-\eta)d_w(\eta) = G ((\lambda - \nu) \circ F(d_w(\eta))),
\]
where $\circ$ means the Hadamard product, $\ie$ pointwise product of two vectors, and $F$ is the cdf of renewable power.

When $K < \min_i d_i(\eta)$, we can see that $d(\eta)$ satisfies the optimal condition therefore $d_w(\eta) = d(\eta)$, $\Delta \overline{\textsf{\text{cs}}}(\eta)=0$. Otherwise, $d_w(\eta) = d(\eta) + G\delta/(2-\eta)$, where $\delta = (\lambda - \nu) \circ (1 - F(d_w(\eta) )\ge 0$.
\[
\begin{array}{rl}
\Delta \overline{\textsf{\text{cs}}}(\eta) & = \frac{1}{2}\{(d_w(\eta))^{\T} G^{-1} d_w(\eta) - (d(\eta))^{\T} G^{-1} d(\eta)\} \\
& = \frac{1}{2}\{ 2 \delta^{\T} d(\eta) +  \delta^{\T} G \delta \} > 0. \\
\end{array}
\]
For the RP, for all $K$ and $\eta$,
\[
\begin{array}{rl}
\Delta \overline{\textsf{\text{rp}}}(\eta)& =(1-\eta)\{(d_w(\eta))^{\T} G^{-1} d_w(\eta) - (d(\eta))^{\T} G^{-1} d(\eta)\} \\
&~~+ \frac{1}{2K}\{(d_w(\eta))^{\T} \Lambda d_w(\eta)\}>0,
\end{array}
\]
where $\Lambda = \text{diag}(\bar{\lambda}_1,...\bar{\lambda}_{24})$. As $K$ goes to infinity, $d_w(\eta)$ is bounded, $\frac{\Delta \overline{\textsf{\text{cs}}}(\eta)}{\Delta \overline{\textsf{\text{cs}}}(\eta)+\Delta \overline{\textsf{\text{rp}}}(\eta)}$ goes to $\frac{1}{3-2\eta}$.

\begin{IEEEbiography}[{\includegraphics[width=1in,height=1.25in,clip,keepaspectratio]{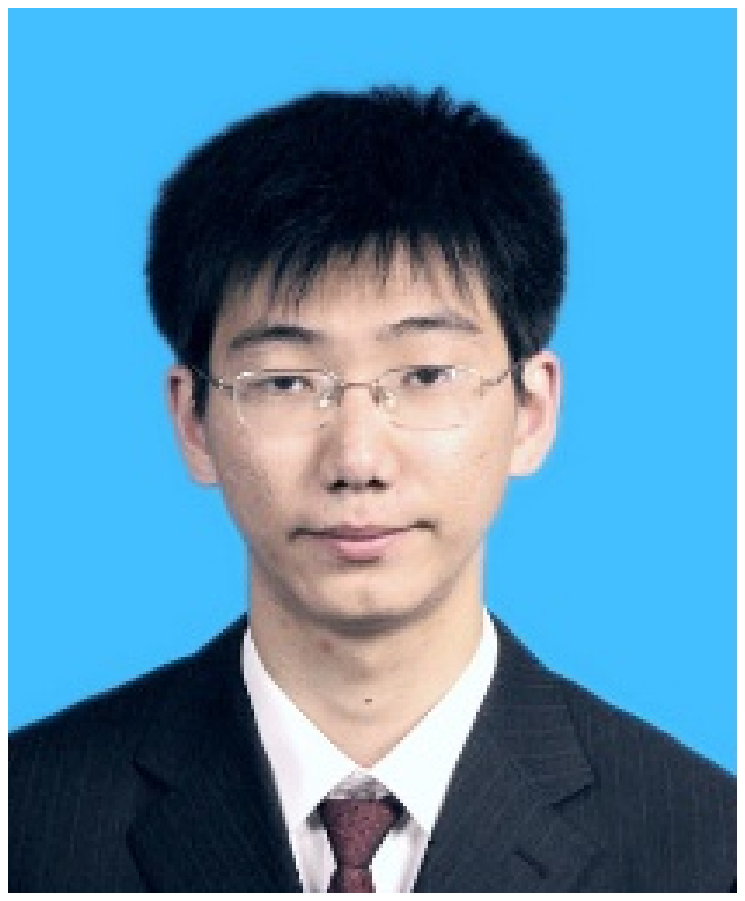}}]{Liyan Jia}
\small
received his B.E. degree from Department of Automation, Tsinghua University in 2009. He is currently working toward the Ph.D. degree in the School of Electrical and Computer Engineering, Cornell University. His current research interests are in smart grid, electricity market and demand response.
\end{IEEEbiography}

\begin{IEEEbiography}[{\includegraphics[width=1in,height=1.25in,clip,keepaspectratio]{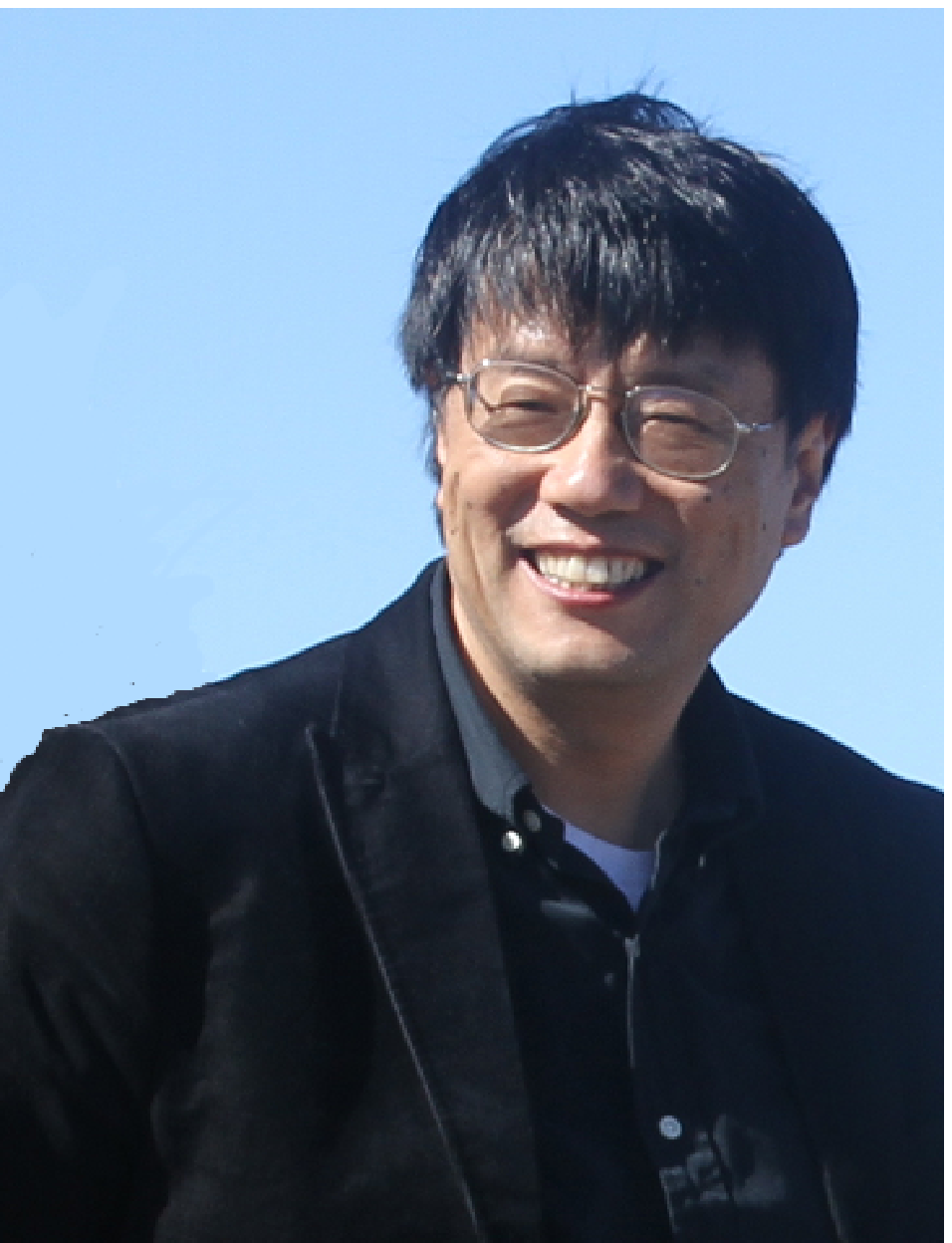}}]{Lang Tong}
\small
\noindent
(S'87,M'91,SM'01,F'05) is the Irwin and Joan
Jacobs Professor in Engineering of Cornell University and the site director of Power Systems Engineering Research Center (PSERC). He received the B.E. degree from Tsinghua University in 1985, and M.S. and Ph.D.
degrees in electrical engineering in 1987 and 1991,
respectively, from the University of Notre Dame.
He was a Postdoctoral Research Affiliate at the Information
Systems Laboratory, Stanford University in 1991.
He was  the 2001 Cor Wit Visiting Professor at
the Delft University of Technology and had held
visiting positions at Stanford University and the University of California at Berkeley.

Lang Tong's research is in the general area of statistical
inference, communications, and complex networks.  His current research focuses on
inference, optimization, and economic problems in energy and power systems.
He received the 1993 Outstanding Young
Author Award from the IEEE Circuits and Systems Society,
the 2004 best paper award  from IEEE Signal Processing Society,
and the 2004 Leonard G. Abraham Prize Paper Award from the
IEEE Communications Society. He is also  a coauthor of seven student paper awards.
He received Young Investigator Award from the Office of Naval Research.
He was a Distinguished Lecturer of the IEEE Signal Processing Society.
\end{IEEEbiography}

\end{document}